\documentclass[oneside,reqno]{amsart}

\usepackage{amssymb}
\usepackage{amsmath}
\usepackage{amsthm}
\usepackage{amsbsy}
\usepackage{bm}
\usepackage{hyperref}
\date{\today}
\usepackage{cite}
\usepackage{array}
\usepackage{xcolor}
\usepackage{tikz}
\usepackage{graphics}
\usepackage{subcaption}
\usepackage{bbm}

\theoremstyle{theorem}
    \newtheorem{theorem}{Theorem}
    \newtheorem{lemma}[theorem]{Lemma}

\theoremstyle{definition} 
    \newtheorem{definition}[theorem]{Definition}
    
    \newtheorem{result}[theorem]{Result}
    \newtheorem{remark}[theorem]{Remark}
    \newtheorem{example}[theorem]{Example}
    \newtheorem{exercise}[theorem]{Exercise}

\def\suchthat{\; : \;}

\def\C{\mathbb{C}}

\def\N{\mathbb{N}}

\def\v{{\bf v}}

\def\<{\langle}
\def\>{\rangle}

\newcommand{\E}{\mbox{E}}

\def\bar{\overline}

\newcommand\Tr{{\mbox{Tr}}}

\newcommand\mnote[1]{} 
\newcommand\be{\begin{equation*}}

\newcommand\ee{\end{equation*}}

\newcommand\ben{\begin{equation}}
\newcommand\een{\end{equation}}
\newcommand\bes{\begin{eqnarray*}}
\newcommand\ees{\end{eqnarray*}}

\newcommand\bex{\begin{exercise}}
\newcommand\eex{\end{exercise}}
\newcommand\beg{\begin{example}}
\newcommand\eeg{\end{example}}
\newcommand\benu{\begin{enumerate}}
\newcommand\eenu{\end{enumerate}}
\newcommand\beit{\begin{itemize}}
\newcommand\eeit{\end{itemize}}
\newcommand\berk{\begin{remark}}
\newcommand\eerk{\end{remark}}
\newcommand\bdefn{\begin{defintion}}
\newcommand\edefn{\end{definition}}
\newcommand\bthm{\begin{theorem}}
\newcommand\ethm{\end{theorem}}
\newcommand\bprf{\begin{proof}}
\newcommand\eprf{\end{proof}}
\newcommand\blem{\begin{lemma}}
\newcommand\elem{\end{lemma}}

\newcommand{\sm}{{\raise0.3ex\hbox{$\scriptstyle \setminus$}}}

\def\CHI{\mathchoice%
{\raise2pt\hbox{$\chi$}}%
{\raise2pt\hbox{$\chi$}}%
{\raise1.3pt\hbox{$\scriptstyle\chi$}}%
{\raise0.8pt\hbox{$\scriptscriptstyle\chi$}}}
\def\smalloplus{\raise1pt\hbox{$\,\scriptstyle \oplus\;$}}

\title[Limiting spectral distribution and joint convergence]{On limiting spectral distribution and joint convergence of some patterned random matrices}

\author{Shambhu Nath Maurya \\\\
Department of Mathematics, 
Indian Institute of Technology Bombay\\
Powai, Mumbai 400076, India}

\date{\today}
\thanks{snmaurya [at] math.iitb.ac.in}

\begin{document}

\begin{abstract}
This article deals with the limiting spectral distribution and joint convergence of reverse circulant and symmetric circulant matrices with independent entries. These results are already proved in articles Bose and Sen (2008) \cite{bose_sen_LSD_EJP}, and Bose, Hazra and Saha (2011) \cite{bose_saha_patter_JC_annals}, respectively. But this article provides alternate proofs of these results which are shorter than the existing proofs. Our method is mainly based on moment method.
\end{abstract}

\maketitle

\noindent{\bf Keywords :} Reverse circulant matrix, symmetric circulant matrix, limiting spectral distribution, expected spectral distribution, joint convergence, moment method, trace-moment formula.
\vskip5pt
\noindent{\bf AMS 2010 subject classification:} 60B20, 60B10, 60C05.

\section{Introduction and main results}
Let $A_n$ be an $n\times n$ real symmetric random matrix. Then \textit{empirical spectral distribution} (ESD) of $A_n$ is defined as
\begin{equation} \label{eqn:ESD}
F_{A_n} (x) = \frac{1}{n}\sum_{k=1}^{n} \mathbbm{1} (\lambda_k \leq x),
\end{equation}
where $\lambda_1,\lambda_2,\ldots, \lambda_n$ are eigenvalues of $A_n$. The weak limit of ESD is known as \textit{limiting spectral distribution} (LSD) of $A_n$, if it exists, in probability or almost surely.

The LSD of random matrices is a quite interesting area of research in Random matrix theory. 
There are several existing methods to study the LSD of random matrices.
One of the crucial methods is the method of moments (see \cite{bryc_Lsd_anna}, \cite{bose_sen_LSD_EJP}) which deals with mainly symmetric matrices. 
 On the other hand, to deal with non-symmetric matrices, Pastur’s fundamental technique of Stieltjes transforms (see \cite{pastur_Pcondtion_first}, \cite{Schatterjee_pastur_stijes}) and the theory of free probability (see \cite{mingo_freeprob_book}, \cite{bose_patterned}) are the important methods.
For various results on LSD of Wigner, sample covariance matrices and other patterned matrices,  we refer the readers to \cite{pastur_LSDfirst_Smatrix}, \cite{bai_jack_LMRA_book}, \cite{bose_patterned}, \cite{bosesaha_circulantbook}. 

 In 1991, Voiculescu \cite{Dan_wigner_JC} studied the joint convergence of Wigner matrices and introduced the notion of free independence. After that many aspects of joint convergence of various random matrices have been studied, see \cite{bose_saha_patter_JC_annals}, \cite{mingo_freeprob_book}, \cite{bose_patterned}.

This article is focused on the LSD and joint convergence of reverse circulant and symmetric circulant matrices whose entries are independent and satisfy some moment condition.
Reverse circulant and symmetric circulant matrices are defined as follows.

\noindent \textbf{Reverse circulant matrix:} An $n\times n$ reverse circulant matrix is defined as
$$
RC_n=\left(\begin{array}{cccccc}
x_1 & x_2 & x_3 & \cdots & x_{n-1} & x_n \\
x_2 & x_3 & x_4 & \cdots & x_{n} & x_{1}\\
x_3 & x_4 & x_5 & \cdots & x_{1} & x_{2} \\
\vdots & \vdots & {\vdots} & \ddots & {\vdots} & \vdots \\
x_n & x_1 & x_2 & \cdots & x_{n-2} & x_{n-1}
\end{array}\right).
$$
 For $j=1,2,\ldots, (n-1)$, its $(j+1)$-th row is obtained by giving its $j$-th row a left circular shift by one position. Note that $RC_n$ is a symmetric matrix and its (i,\;j)-th element is $x_{(i+j-1) \bmod  n}$. 
\vskip5pt
\noindent\textbf{Symmetric circulant matrix:} An $n\times n$  symmetric circulant matrix is defined as 
$$
SC_n=\left(\begin{array}{cccccc}
x_0 & x_1 & x_2 & \cdots & x_{2} & x_1 \\
x_1 & x_0 & x_1 & \cdots & x_{3} & x_{2}\\
x_{2} & x_1 & x_0 & \cdots & x_{4} & x_{3}\\
\vdots & \vdots & {\vdots} & \ddots & {\vdots} & \vdots \\
x_1 & x_2 & x_3 & \cdots & x_1 & x_0
\end{array}\right).
$$
  For $j=1,2,\ldots, (n-1)$, its $(j+1)$-th row is obtained by giving its $j$-th row a right circular shift by one position. Note that $SC_n$ is a symmetric matrix and its (i,\;j)-th element is $x_{\frac{n}{2}-|\frac{n}{2}-|i-j||}$.

The following results provide the LSD of $RC_n$ and $SC_n$, respectively.
\begin{theorem}\label{thm:revESD}
 Suppose $RC_n$ is the reverse circulant matrix with entries $\{\frac{X_j}{\sqrt n};{j\geq 0}\}$, where $\{X_j\}$ are i.i.d. with mean zero and variance one. Then LSD of $RC_n$ is  \textit{symmetrized Rayleigh distribution} whose density function is $ |x| \exp(-x^2)$ for   $\ -\infty < x < \infty.$
\end{theorem}
\begin{theorem}\label{thm:symESD}
 Suppose  $SC_n$ is the symmetric circulant matrix with entries $\{\frac{X_j}{\sqrt n};{j\geq 0}\}$, where $\{X_j\}$ are i.i.d. with mean zero and variance one.  Then LSD of $SC_n$ is the standard Gaussian distribution.
\end{theorem}
 Theorem \ref{thm:revESD} and Theorem \ref{thm:symESD} are proved in \cite{bose_sen_LSD_EJP} using moment method. We shall also use moment method but our combinatorial techniques are different than used in \cite{bose_sen_LSD_EJP}. Note that the $h$-th moment of ESD of a real symmetric matrix $A_n$ can be written in terms of trace of $(A_n)^h$ by using \textit{trace-moment formula}. Thus, if closed form of trace of $(A_n)^h$ is known then one can exploit it to study the LSD of $A_n$. We derive a convenient trace formulas for $(RC_n)^h$ and $(SC_n)^h$, and then use some nice combinatorial arguments to establish the results.  
 
Now we state an assumption for a sequence of random variables $\{X_j\}_{j\geq 0}$ which will be required for the joint convergence for  reverse circulant and symmetric circulant matrices.
\begin{equation}\label{eqn:condition_ESD}
\{X_j\} \mbox{ are independent}, \ \E(X_j)=0, \ \E(X_j^2)=1 \ \forall \ j \geq 0 \mbox{ and }  \sup_{j\geq 0}\E(|X_j|^k)=\alpha_k<\infty \ \mbox{for}\ k\geq 3.
\end{equation}
The following results provide the joint convergence of reverse circulant and symmetric circulant matrices, respectively.
\begin{theorem}\label{thm:rev_JC}
Suppose  $\{RC_n^{(i)} \}_{ 1 \leq i \leq m}$ are $m$ independent reverse circulant matrices whose entries are $\{\frac{X_j^{(i)} }{\sqrt n};{j\geq 0}\}_{ 1 \leq i \leq m}$, where for $i=1, 2, \ldots, m$, $\{X_j^{(i)}\}_{j\geq 0}$ satisfy assumption \eqref{eqn:condition_ESD}.   
 Then as elements of a non-commutative probability space $(\mathcal{M}_n, \phi_n = \frac{1}{n} \E \Tr )$,
 \begin{equation*} \label{eqn:joint_conv_rev}
\{RC_n^{(i)} \}_{ 1 \leq i \leq m} \stackrel{J}{\rightarrow} \{a_i\}_ {1 \leq i \leq m},
 \end{equation*}
 where $(\mathcal{M}_n, \phi_n = \frac{1}{n} \E \Tr )$ and $\{a_i\}_ {1 \leq i \leq m}$ are as given in Definition \ref{def:non-commu} and Example \ref{exam:example_rJC}, respectively. Here $\stackrel{J}{\rightarrow}$ denotes the joint convergence of random matrices as given in Definition \ref{def:joint_convergence}.
 \end{theorem} 
\begin{theorem}\label{thm:sym_JC}
Suppose $\{SC_n^{(i)} \}_{ 1 \leq i \leq m}$ are $m$ independent symmetric circulant matrices whose entries are $\{\frac{X_j^{(i)} }{\sqrt n};{j\geq 0}\}_{ 1 \leq i \leq m}$, where for $i=1, 2, \ldots, m$, $\{X_j^{(i)}\}_{j\geq 0}$ satisfy assumption \eqref{eqn:condition_ESD}.   
 Then as elements of the non-commutative probability space $(\mathcal{M}_n, \phi_n = \frac{1}{n} \E \Tr )$,
 \begin{equation*} \label{eqn:joint_conv_sym}
\{SC_n^{(i)} \}_{ 1 \leq i \leq m} \stackrel{J}{\rightarrow} \{b_i\}_ {1 \leq i \leq m},
 \end{equation*}
where $\{b_i\}_{1 \leq i \leq m}$ are independent standard Gaussian variables.
 \end{theorem} 
 Theorem \ref{thm:rev_JC} and  Theorem \ref{thm:sym_JC} are proved in \cite{bose_saha_patter_JC_annals} but here we provide alternate proofs, by using different combinatorial techniques and convenient trace formulas for product of these matrices.

Here is a brief outline of the rest of the manuscript. In Section \ref{pre_ESD_RcSc} we state some basic notation and a result which are required for the proofs of the theorems. In Section \ref{sec:ESD_revsym} we prove Theorem \ref{thm:revESD} and Theorem \ref{thm:symESD} by using moment method and some nice combinatorial arguments. Finally in Section \ref{sec:JC_revsym} we complete the proofs of Theorem \ref{thm:rev_JC} and Theorem \ref{thm:sym_JC}.

 \section{Preliminaries}\label{pre_ESD_RcSc}
Here we recall some standard notation and state a result.  Suppose $F_{A_n}$ is the ESD of a real symmetric matrix $A_n$. Then
\begin{equation}\label{eqn:moment_Tr_formula}
\tilde{E}_h(F_{A_n}) = \frac{1}{n} \sum_{k=1}^n (\lambda_k)^h = \frac{1}{n} \Tr(A_n)^h,
\end{equation}     
      where $\tilde{E}_h$ denotes the $h$-th moment with respect to the \textit{empirical spectral measure} of $A_n$ 
       and $\Tr$ denotes the trace of a matrix. The above expression is usually known as \textit{trace-moment formula}.
      
Now we state a result from Chapter 2 of \cite{bosesaha_circulantbook} which provides sufficient condition for the convergence of ESD 
of a given real symmetric random matrix $A_n$.
\begin{result} \textbf{(Lemma 2.2.4, \cite{bosesaha_circulantbook})} \label{res:M_1234} 
 Suppose there exists a sequence $\{\beta_h\}_{h \geq 1}$ such that 
\begin{enumerate}
		\item [$(M_1)$] for every $h\in \mathbb{N}$, $\E[ \tilde{E}_h(F_{A_n})] \rightarrow \beta_h$ as $n \to \infty$,
		\item [$(M_2)$] there is a unique distribution $F$ whose moment sequence is $\{\beta_h\}_{h \geq 1}$, and
		\item [$(M_3)$] for every $h\in \mathbb{N}$, $\sum_{n=1}^{\infty}\E \big[ \tilde{E}_h(F_{A_n})-\E[ \tilde{E}_h(F_{A_n})] \big]^4<\infty$.
	\end{enumerate}
	Then weakly, $\{F_{A_n}\}$ converges almost surely to $F$.
\end{result}
    Now we recall the notion of non-commutative probability space and joint convergence. 
    \begin{definition} \label{def:non-commu}
    Let $\mathcal{A}$ over $\C$ be a unital algebra  with unity $1_A$. Suppose $\phi$ is a linear functional from $\mathcal{A}$ to $\C$ with $\phi(1_A) = 1$. Then $(\mathcal{A}, \phi)$ is called a \textit{non-commutative probability space} and $\phi$ is called a \textit{state}.

    Suppose $\mathcal{M}_n$ is the set of all $n\times n$ random matrices whose entries have finite moments of all order and the state $\phi_n$ is defined as
$$ \phi_n(A_n) = \frac{1}{n} \E[\Tr(A_n)]  \ \ \mbox{ for } A_n \in \mathcal{M}_n.$$
Then $(\mathcal{M}_n, \phi_n)$ is a non-commutative probability space.
    \end{definition}
       \begin{definition}\label{def:joint_convergence}
Let $(\mathcal{A}_n, \tau_n), \ n \geq 1$ and $(\mathcal{A}, \tau)$ be non-commutative probability space. Suppose $\{a^{(i)}_n\}_{i\in I} \subset \mathcal{A}_n$, where $I$ is any finite subset of $\N$. Then we say $\{a^{(i)}_n\}_{i\in I}$ \textit{jointly converges} to $\{a^{(i)}\}_{i\in I} \subset \mathcal{A}$, denoted by $\{a^{(i)}_n\}_{i\in I}  \stackrel{J}{\rightarrow} \{a^{(i)}\}_{i\in I}$, if for every complex polynomial $Q$, 
	$$\tau_n(Q(a^{(i)}_n; i \in I)) \to \tau(Q(a^{(i)}; i \in I)) \ \mbox{  as $n \to \infty$}.$$     
 \end{definition}
    
\section{Proofs of Theorem \ref{thm:revESD}  and Theorem \ref{thm:symESD}} \label{sec:ESD_revsym}
First we state a result which provides an alternate condition on the entries of $RC_n$ and $SC_n$, under which LSD of $RC_n$ and $SC_n$ exist.
\begin{result}   \textbf{(Theorem 2, \cite{bose_sen_LSD_EJP})}. \label{res:iid_uniform_ESD}
Let $A_n$ be a reverse circulant or symmetric circulant matrix with entries  $\{\frac{X_j}{\sqrt n};{j\geq 0}\}$. Suppose ESD of $A_n$ converges almost surely to some fixed non random distribution $F$ when $\{X_j\}$ satisfy assumption (\ref{eqn:condition_ESD}).  Then the same limit continues to hold if $\{X_j\}$ are i.i.d. with mean zero and variance one.
\end{result}
Now we state a  convenient trace formula for reverse circulant matrices which will be required to prove Theorem \ref{thm:revESD} and again in Section \ref{sec:JC_revsym} to prove Theorem \ref{thm:rev_JC}.
\begin{lemma} \label{lem:traceRc_copy}
 Suppose $\{RC_n^{(r)}\}_{1\leq r \leq h}$ are the reverse circulant matrices with entries $\{\frac{X_j^{(r)} }{\sqrt n};{j\geq 0}\}$. Then 
\begin{equation}\label{trace formula RC_n}
\Tr[ RC_n^{(1)} RC_n^{(2)} \cdots RC_n^{(h)}]
=  \left\{\begin{array}{ll} 
		 	 \displaystyle \frac{1}{n^{p-1}} \sum_{A_{2p}} X_{j_1}^{(1)} X_{j_2}^{(2)} \cdots X_{j_{2p}}^{(2p)}  & \text{if }       h=2p, \\
		 	 \displaystyle  \frac{1}{n^{ \frac{2p+1}{2} }} \sum_{i=1}^n \sum_{A_{p, i}} X_{j_1}^{(1)}   X_{j_2}^{(2)} \cdots X_{j_{2p+1}} ^{(2h+1)} & \text{if } h=2p+1,	 
		 	  \end{array}\right. 
 \end{equation}  
where 
\begin{align}\label{def:A_2p_rESD}
A_{2p} &=\big\{(j_1, j_2, \ldots, j_{2p})\in \mathbb N^{2p}\suchthat \sum_{k=1}^{2p}(-1)^k j_k=0 \mbox{ (mod $n$)}, 1\le j_1,\ldots, j_{2p}\le n\big\}, \\
A_{p,i} &=\big\{(j_1, j_2, \ldots, j_{2p+1})\in \mathbb N^{2p+1}\suchthat \sum_{k=1}^{2p+1}(-1)^k j_k= 2i-1 \mbox{ (mod $n$)}, 1\le j_1,\ldots, j_{2p+1}\le n\big\}. \nonumber
\end{align}
\end{lemma}
For the proof of above lemma, we refer the reader to \cite{m&s_RcSc_jmp}. 
Now we recall a notion of {\it odd-even pair matched} vector from \cite{m&s_RcSc_jmp}.
\begin{definition}\label{def:odd-even_match_rESD}
	Suppose $v=(j_1, j_2, \ldots, j_p)$ is a vector in $\mathbb N^p$. Two elements $j_k, j_\ell$ of $v$ are said to be {\it odd-even pair matched} if $j_k=j_\ell$, $j_k$ appears exactly twice in $v$, once each at an odd and an even position. For example, in $(5,5,3,4)$, entry $5$ is {\it odd-even pair matched} whereas in $(5,2,5,3)$, entry $5$ is  not {\it odd-even pair matched}.
	We shall call $v$ {\it odd-even pair matched} if all its entries are {\it odd-even pair matched}. In that case, $p$ is necessarily even.
\end{definition}
Note that, if $J_{2p}=(j_1,j_2,\ldots, j_{2p})\in A_{2p}$  and the entries $\{j_1,j_2,\ldots, j_{2p} \}$ are at least pair matched, then the maximum number of free entries in $J_{2p}$ will be $p$, only when $J_{2p}$ is {\it odd-even pair matched}. By free entries we mean that these entries can be chosen freely from their range. 
This observation will be used in the proofs of Theorem \ref{thm:revESD} and Theorem \ref{thm:rev_JC}. 
 
Now with the help of above result and trace formula, we prove Theorem \ref{thm:revESD}.
\begin{proof}[Proof of Theorem \ref{thm:revESD}] In view of Result \ref{res:iid_uniform_ESD}, it is enough to assume (\ref{eqn:condition_ESD}) for the entries $\{X_j\}$ of $RC_n$ and prove the required almost sure convergence.
 First recall from Result \ref{res:M_1234} that to prove Theorem \ref{thm:revESD}, it is sufficient to verify $(M_1), (M_2)$ and $(M_3)$ for $RC_n$. 

First we verify $(M_1)$.
Note from (\ref{trace formula RC_n}) that for odd and even values of $h$ we have different trace formulae for $(RC_n)^h$. So we calculate the limit of $\E[ \tilde{E}_h(F_{RC_n})]$ for odd and even values of $h$ separately.
First suppose $h$ is even, say $h=2p$. Then from (\ref{eqn:moment_Tr_formula}) and (\ref{trace formula RC_n}), we get
\begin{align} \label{eqn:E_X_i_rev}
\E\big[ \tilde{E}_{2p}(F_{RC_n}) \big] = \frac{1}{n} \E[ \Tr(RC_n)^{2p}]=  \frac{1}{n^p} \sum_{A_{2p}} \E[ X_{j_1} X_{j_2} \cdots X_{j_{2p}}],
\end{align} 
where $(j_1, j_2, \ldots, j_{2p}) \in A_{2p}$ and $A_{2p}$ is as defined in (\ref{def:A_2p_rESD}).

 Observe from the assumption (\ref{eqn:condition_ESD}) that $\E[X_j] =0$ for each $j$, therefore the right side of (\ref{eqn:E_X_i_rev}) will be non-zero when the entries $\{ j_1, j_2, \ldots, j_{2p}\}$ are at least pair matched. Now with assumption (2), if we combine the above observation with Definition \ref{def:odd-even_match_rESD}, we get
 \begin{align}  \label{eqn:E_bound_rEven}
\sum_{A_{2p}} \E[ X_{j_1} X_{j_2} \cdots X_{j_{2p}}]
&= \left\{\begin{array}{ll} 
		 	O(n^p) & \text{if} \    (j_1, j_2, \ldots, j_{2p}) \mbox{ is odd-even pair matched}, \\
		 	 o(n^p) & \text{otherwise}. 	 
		 	  \end{array}\right. 
\end{align}
Hence from (\ref{eqn:E_X_i_rev}), we get
\begin{align*}
\lim_{n\to\infty} \E\big[ \tilde{E}_{2p}(F_{RC_n}) \big] = \lim_{n\to\infty} \frac{1}{n^p}  (n^p p!) = p!, 
\end{align*}
 where the first equality arises due to the odd-even pair matching among the entries of $(j_1, j_2, \ldots, j_{2p})$. 
 Out of $2p$ many entries, if we freely choose $p$ entries in odd positions, then the rest $p$ many entries in even positions can permute among themselves in $p!$ ways. Hence odd-even pair matching among the entries of $(j_1, j_2, \ldots, j_{2p})$ happen in $n^pp!$ ways.
 
 Now suppose $h$ is odd, say $h=2p+1$. Then from the trace formula (\ref{trace formula RC_n}), we get
\begin{align} \label{eqn:E_X_i_rev_odd}
\E\big[ \tilde{E}_{2p+1}(F_{RC_n}) \big] & = \frac{1}{n^\frac{2p+3}{2}} \sum_{i=1}^n \sum_{A_{p, i}} \E[ X_{j_1} X_{j_2} \cdots X_{j_{2p+1}}],
\end{align} 
where for each $i=1, 2, \ldots, n$, $(j_1, j_2, \ldots, j_{2p+1}) \in A_{p,i}$ and $A_{p,i}$ is as defined in (\ref{def:A_2p_rESD}).

Since $\E[X_j] =0$, $\E[ X_{j_1} X_{j_2} \cdots X_{j_{2p+1}}]$ will be non-zero when the entries $\{ j_1, j_2, \ldots, j_{2p+1}\}$ are at least pair matched. 
Thus a typical term of (\ref{eqn:E_X_i_rev_odd}) has the maximum contribution if one entry is triple matched and the rest entries are pair-matched. In this case, the contribution will be of the order $O(n^{\frac{2p}{2} -1})$, 
 where $(-1)$ arises due to the constraint of $A_{p,i}$, $\sum_{k=1}^{2p+1} (-1)^k j_k=2i-1 \mbox{ (mod $n$)}$. 
So, we get
 \begin{align} \label{eqn:ESDo(1)_rESD} 
\E\big[ \tilde{E}_{2p+1}(F_{RC_n}) \big] 
= \frac{1}{n^\frac{2p+3}{2}} O(n)  O(n^{p-1})= o(1).
\end{align}
Hence
  \begin{align*}
\lim_{n\to\infty}  \E\big[ \tilde{E}_{h}(F_{RC_n}) \big] 
 = 
 \left\{\begin{array}{ll} 
		 	 p! & \text{if }       h=2p, \\
		 	 0 & \text{if } h=2p+1. 	 
		 	  \end{array}\right. 
	 \end{align*}
	 
	 Observe from the above limit that the limiting sequence satisfy the condition $(M_2)$ and it corresponds to a unique distribution which is known as \textit{symmetrized Rayleigh distribution}.

Now we check condition $(M_3)$.
First note that to verify $(M_3)$, it is sufficient to show that for every $h \in \mathbb{N}$,
\begin{equation} \label{eqn:M_3_r}
\E \big[ \tilde{E}_h(F_{RC_n})-\E[ \tilde{E}_h(F_{RC_n})] \big]^4 = O(n^{-2}).
\end{equation}
Now we define a term called {\it connected vectors} which will be required to establish (\ref{eqn:M_3_r}).
\begin{definition}\label{def:conn_graph}
Suppose $T=\{J^1, J^2, \ldots, J^\ell\}$ is a set of $\ell$ vectors,  where $J^i =(j^i_1, j^i_2, \ldots, j^i_{p})\in \mathbb N^{p}$ for $1 \leq i \leq \ell$. Then we can associate a graph from the given vectors in following scheme: $J^i$ denotes the vertex and there will be an edge between $J^r$ and $J^s$ if  $ \{\{j^r_1, j^r_2, \ldots, j^r_{p}\} \cap \{j^s_1, j^s_2, \ldots, j^s_{p}\}  \}  \neq \emptyset$.
A subset $T'=\{J^{i_1}, J^{i_2}, \ldots, J^{i_d}\}$ of $T$ is said to be \textit{connected} in $T$ if the graph formed by $T'$ is a connected component of the graph formed by $T$.
\end{definition}
As before, we derive (\ref{eqn:M_3_r}) for odd and even values of $h$ separately. 
First suppose $h$ is even, say $h=2p$. Then from the trace formula (\ref{trace formula RC_n}), we get
\begin{align} \label{eqn:E^4_rEven}
\E \big[ \tilde{E}_{2p} (F_{RC_n})-\E[ \tilde{E}_{2p} (F_{RC_n})] \big]^4 = \frac{1}{n^{4p} } \sum_{A_{2p}, A_{2p}, A_{2p}, A_{2p} }  \E \Big[ \prod_{r=1}^4 \big( X_{J^r_p}  -  \E[X_{J^r_p}] \big) \Big] = \frac{1}{n^{4p} } L, \mbox{ say},
\end{align}
where for each $r=1, 2, 3, 4$, $J^r_{p} = (j^r_1, j^r_2, \ldots, j^r_{2p}) \in A_{2p}, A_{2p}$ as in (\ref{def:A_2p_rESD}) and $X_{J^r_p} = X_{j^r_1} X_{j^r_2} \cdots X_{j^r_{2p}} $. 

Now we calculate the order of convergence of $L$. Depending on the connectedness 
between $J^r_{p}$'s, the following three cases arise:
\vskip5pt
\noindent \textbf{Case I.} \textbf{ At least one $J^r_{p}$ is not connected with the remaining one:}
 In this case, we get $L=0$ due to the independence of the entries $\{X_j\}$.
\vskip5pt
\noindent \textbf{Case II.} \textbf{$J^1_{p}$ is connected with one of $J^2_{p}, J^3_{p}, J^4_{p}$ only and the remaining two of $J^2_{p}, J^3_{p}, J^4_{p}$  are also connected with themselves only:} Without loss of generality, we assume $J^1_{p}$ is connected with $J^2_{p}$ only and $J^3_{p}$ is connected with $J^4_{p}$ only.
 So from the independence of the entries, we get
\begin{align} \label{eqn:case1_Sr}
L =  \sum_{ A_{2p}, A_{2p}} \E \big[ \prod_{r=1}^2 \big( X_{J^r_p}  -  \E[X_{J^r_p}] \big) \big]  \sum_{A_{2p}, A_{2p}} \E \big[\prod_{r=3}^4 \big( X_{J^r_p}  -  \E[X_{J^r_p}] \big) \big].
\end{align}
Now under assumption (\ref{eqn:condition_ESD}), there exists a positive constant $\alpha$ such that  
 \begin{align*} 
\Big| \sum_{A_{2p}, A_{2p}}   \E\big[ \prod_{r=1}^2 \big( X_{J^r_p}  -  \E[X_{J^r_p}] \big) \big] \Big| = \Big| \sum_{( J^1_{p}, J^2_{p}) \in B_{P_2}}  \E\big[ \prod_{r=1}^2 \big( X_{J^r_p}  -  \E[X_{J^r_p}] \big) \big] \Big| \leq \alpha |B_{P_2}|,
\end{align*}
where  $B_{P_2} \subseteq A_{2p} \times A_{2p}$ such that $( J^1_{p}, J^2_{p}) \in B_{P_2}$ if $J^1_{p}$ is connected with $J^2_{p}$ and
the entries of $\{\{j^1_1,\ldots, j^1_{2p}\} \cup \{ j^2_1, \ldots, j^2_{2p}\} \}$ are at least pair matched. Note that the above inequality arises due to  uniform boundedness of moments of $X_j$'s. Here $|\{\cdot\}|$ denotes the cardinality of the set $\{\cdot\}$.

 It is easy to observe from Definition \ref{def:odd-even_match_rESD} that the maximum $(2p-1)$ many entries can be chosen freely in $B_{P_2}$. Hence $|B_{P_2}|= O(n^{2p-1}).$ 
Now with this observation, if we combine the above inequality with (\ref{eqn:case1_Sr}), we get $|L| \leq {\alpha}^2  \big( |B_{P_2}| \big)^2 = O(n^{4p-2})$.
\vskip5pt
 \noindent \textbf{Case III.} \textbf{$J^1_{p}, J^2_{p}, J^3_{p}, J^4_{p}$ are connected:}  
By the argument similar as used in Case II, we get
\begin{align*} 
|L| =  \Big| \sum_{( J^1_{p}, J^2_{p}, J^3_{p}, J^4_{p}) \in B_{P_4}} \E \Big[ \prod_{r=1}^4 \big( X_{J^r_p}  -  \E[X_{J^r_p}] \big) \Big]  \Big| \leq \beta |B_{P_4} |,
\end{align*}
 where $\beta$ is a positive constant and $B_{P_4}$ is defined as
 \begin{align*} 
 B_{P_4}=\big\{( J^1_{p}, J^2_{p}, J^3_{p}, J^4_{p}) \in  &  \ A_{2p} \times A_{2p} \times A_{2p} \times A_{2p}  :  J^1_{p}, J^2_{p}, J^3_{p}, J^4_{p} \mbox{ are connected and}  \\
 & \ \qquad \mbox{the entries of }  \cup_{r=1}^4 \{ j^r_1, j^r_2, \ldots, j^r_{2p}\} \mbox{ are at least pair matched} \big\}.
\end{align*}
 Note that $|B_{P_4}| = o({n}^{4p-2})$. For the details about the cardinality of $B_{P_4}$, see \cite[Lemma 18]{m&s_RcSc_jmp}. Therefore from the above inequality, we get $|L| \leq o(n^{4p-2})$.

Hence on combining all three cases, we get 
\begin{equation} \label{eqn:E^4_O(4p-2)_r}
|L| \leq  \sum_{A_{2p}, A_{2p}, A_{2p}, A_{2p} }  \Big|  \E \Big[ \prod_{r=1}^4 \big( X_{J^r_p}  -  \E[X_{J^r_p}] \big) \Big]  \Big| 
 \leq O(n^{4p-2}).
 \end{equation}
Finally, on combining (\ref{eqn:E^4_rEven}) and (\ref{eqn:E^4_O(4p-2)_r}), we get
$$\E \big[ \tilde{E}_{2p} (F_{RC_n})-\E[ \tilde{E}_{2p}(F_{RC_n})] \big]^4 = O(n^{-2}).$$

Now suppose $h$ is odd, say $h=2p+1$. Then from the trace formula (\ref{trace formula RC_n}), we get
\begin{align} \label{eqn:E^4_rOdd}
&\E \big[ \tilde{E}_{2p+1} (F_{RC_n})-\E[ \tilde{E}_{2p+1} (F_{RC_n})] \big]^4 \nonumber \\
& \qquad = \frac{1}{n^{4p+6} } \sum_{i_1, i_2, i_3, i_4=1}^n  \sum_{A_{p, i_1}, A_{p, i_2}, A_{p, i_3}, A_{p, i_4}}  \E \Big[ \prod_{r=1}^4 \big( X_{J_{p,i_r} }  -  \E[X_{J_{p, i_r} }] \big) \Big],
\end{align}
where $J_{p, i_r} = (j^r_1, \ldots, j^r_{2p+1}) \in A_{p, i_r}, A_{p, i_r}$ as in (\ref{def:A_2p_rESD}) and $X_{J_{p, i_r} } = X_{j^r_1} \cdots X_{j^r_{2p+1}}$  for $r=1, 2, 3, 4$. 

If we follow the arguments similar to those used to establish (\ref{eqn:E^4_O(4p-2)_r}), we get
\begin{equation*} 
\sum_{A_{p, i_1}, A_{p, i_2}, A_{p, i_3}, A_{p, i_4}}   \Big|  \E \Big[ \prod_{r=1}^4 \big( X_{J_{p,i_r} }  -  \E[X_{J_{p, i_r} }] \big) \Big] \Big| 
 \leq O(n^{4p}).
 \end{equation*}
Now on combining the above inequality with (\ref{eqn:E^4_rOdd}), we get (\ref{eqn:M_3_r}) for odd values of $h$.
This completes the proof of Theorem \ref{thm:revESD}.
\end{proof}
Now we prove Theorem \ref{thm:symESD}. The idea of the proof is similar to the proof of Theorem \ref{thm:revESD}.
First we recall a convenient trace formula from \cite{m&s_RcSc_jmp} for $SC_n$ with entries $\{\frac{X_j}{\sqrt{n} }; j\geq 0\}$:
\begin{equation}\label{trace,formula_sp}
\Tr(SC_n)^p
   = \left\{\begin{array}{ccc} 	 
		 	\displaystyle \frac{1}{n^{\frac{p-2}{2} }} \sum_{k=0}^{p}\binom{p}{k}X_0^{p-k}\sum_{ C_{k}} X_{J_{k}}  & \text{if}& \mbox{ $n$ is odd},\\
			\displaystyle \frac{1}{2n^{\frac{p-2}{2} }}\sum_{k=0}^{p}\binom{p}{k} \Big[ Y_k  \sum_{ C_k} X_{J_{k}} + \tilde{Y}_k \sum_{ \tilde{C}_k}  X_{J_{k}} \Big] & \text{if}& \mbox{ $n$ is even}, 	 	 
		 	  \end{array}\right.	
 \end{equation}
where  for each  $k=0, 1, 2,\ldots, p$, 
\begin{align} \label{def:A_2p_sESD}
C_{k}&=\Big\{(j_1,\ldots,j_{k})\in\mathbb N^k \hskip-5pt \suchthat \hskip-2pt \sum_{i=1}^{k}\epsilon_i j_i=0\; \mbox{(mod n)}, \epsilon_i\in\{+1,-1\}, 1\le j_1,\ldots,j_{k}\le \Big\lfloor\frac{n}{2} \Big\rfloor \Big\}, \\
\tilde{C}_{k} &=\Big\{(j_1,\ldots,j_{k})\in\mathbb N^k \hskip-5pt \suchthat \hskip-2pt \sum_{i=1}^{k}\epsilon_i j_i=0\; \mbox{(mod } \frac{n}{2}), \ \sum_{i=1}^{k}\epsilon_i j_i \neq 0\; \mbox{(mod }n), \epsilon_i\in\{+1,-1\}, 1\le j_1,\ldots,j_{k} \leq (\frac{n}{2}-1)\Big\},\nonumber \\
Y_k & = {(X_0 +  X_{\frac{n}{2}})}^{p-k} + {(X_0 -  X_{\frac{n}{2}})}^{p-k}, \ \tilde{Y}_k = {(X_0 +  X_{\frac{n}{2}})}^{p-k} - {(X_0 -  X_{\frac{n}{2}})}^{p-k}, 
 X_{J_{k}}  = X_{j_1}\cdots X_{j_{k}}.  \nonumber
\end{align}
where $\lfloor x \rfloor$ denotes the greatest integer less than or equal to $x$.	
Note that $C_0$ and $\tilde{C}_0$ are empty sets with the convention that the contribution from the sum corresponding to $C_0$ and $\tilde{C}_0$ are $1$ in \eqref{trace,formula_sp}. Also in $C_k$ and $\tilde{C}_k$, $(j_1, j_2, \ldots,j_{k})$ are collected according to their multiplicity. 


Now we recall a notion of {\it opposite sign pair matched} vector from \cite{m&s_RcSc_jmp}.
\begin{definition}\label{def:opposite_sign}
Suppose $v=(j_1, j_2, \ldots, j_p)$ is a vector in $\mathbb N^p$ and $(\epsilon_1,\epsilon_2,\ldots,\epsilon_p)\in \{+1,-1\}^p$. Two elements $j_k, j_\ell$ of $v$ are said to be {\it opposite sign pair matched} given $(\epsilon_1,\epsilon_2,\ldots,\epsilon_p)$, if $\epsilon_k$ and $\epsilon_\ell$  are of opposite sign, $j_k= j_\ell$ and $j_k$ appear exactly twice in $v.$
 For example, in $(3,5,8,5)$, $5$ is {\it opposite sign pair matched}, if $\epsilon_2=1$ and $\epsilon_4=-1$ or $\epsilon_2=-1$ and $\epsilon_4=1$ whereas if $\epsilon_2$ = $\epsilon_4= 1$ or $\epsilon_2$ = $\epsilon_4= -1$, then $5$ is not  {\it opposite sign pair matched}. 
For even $p$, we say that a vector $v$ is {\it opposite sign pair matched} for a given $(\epsilon_1,\epsilon_2,\ldots,\epsilon_p)$, if all its entries are {\it opposite sign pair matched}. The set of all such vectors in $C_{p}$ will be denoted by $U_{p}$.	 
\end{definition}
Now we see some observations from Definition \ref{def:opposite_sign} which will be used in the proofs of Theorem \ref{thm:symESD} and Theorem \ref{thm:sym_JC}. For any $p\geq2$, we define a set $C^{(2)}_{p}$ as
\begin{equation} \label{eqn:C^2_k_slsd}
C^{(2)}_{p} = \{ (j_1, j_2, \ldots, j_p) \in C_p  :  \mbox{ the entries of $\{ j_1, j_2, \ldots, j_{p}\}$ are at least pair matched}  \}.
\end{equation}
First note that $ | C^{(2)}_{p}| \leq O(n^{\lfloor \frac{p}{2} \rfloor})$. Now for even values of $p$,  if $(j_1, j_2, \ldots, j_p) \in C^{(2)}_{p}\setminus U_{p}$, then either there exists a $j_i$ which is at least triple matched or there exists a matching pair $j_k,j_\ell$ such that the corresponding $\epsilon_k, \epsilon_\ell$ have same sign. In both the situations, there is a loss of degree of freedom. Hence
\begin{equation} \label{eqn:opposite_sign_card}
 \ \ | (C^{(2)}_{p}\setminus U_{p})| =O(n^{\frac{p}{2}-1})\ \ \text{and} \ \ | U_{p}| =O(n^{\frac{p}{2}}), \ \ \ \ \mbox{whenever $p$ is an even positive integer},
\end{equation}
where $| U_{p}| =O(n^{\frac{p}{2}})$ trivially follows from Definition \ref{def:opposite_sign}. 

Now we prove Theorem \ref{thm:symESD} by using the above observation and the trace formula of $SC_n$. 
\begin{proof}[Proof of Theorem \ref{thm:symESD}]  
The idea of the proof is same as the proof of Theorem \ref{thm:revESD}. We verify the conditions $(M_1), (M_2)$ and $(M_3)$ of Result \ref{res:M_1234} for $SC_n$ under assumption (\ref{eqn:condition_ESD}). 

First we verify $(M_1)$. We calculate the limit of $\E[ \tilde{E}_h(F_{SC_n})]$ for odd and even values of $n$ separately.
First suppose $n$ is odd. Then from (\ref{eqn:moment_Tr_formula}) and (\ref{trace,formula_sp}), we get
\begin{align} \label{eqn:E_X_i_symOdd}
\E\big[ \tilde{E}_{h}(F_{SC_n}) \big] \mathbb I_{\{ n \ odd\}} 
= \frac{1}{n^{\frac{h}{2} }}  \sum_{k=0}^{h}\binom{h}{k} \E[X_0^{h-k}] \sum_{ C_{k}} \E[ X_{j_1} X_{j_2} \cdots X_{j_{k}} ],
\end{align} 
where $J_k =(j_1, j_2, \ldots, j_{k}) \in C_{k}$ and $C_{k}$ is as defined in (\ref{def:A_2p_sESD}).

 Since $\E[X_j] =0$, $\E[ X_{j_1} X_{j_2} \cdots X_{j_{k}} ]$ will be non-zero when $J_k \in C^{(2)}_k$, where $C^{(2)}_k \subseteq C_{k} $ is as defined in (\ref{eqn:C^2_k_slsd}). Also note that for each $k$, the maximum contribution from a typical term of (\ref{eqn:E_X_i_symOdd}) can be of the order $O(n^{ {\lfloor\frac{k}{2}\rfloor} })$ as $ | C^{(2)}_{k}| \leq O(n^{\lfloor \frac{k}{2} \rfloor})$. Therefore (\ref{eqn:E_X_i_symOdd}) will be
 \begin{align} \label{eqn:moment_symOdd}
 \E\big[ \tilde{E}_{h}(F_{SC_n}) \big] \mathbb I_{\{ n \ odd\}}
 = \left\{\begin{array}{ll} 
		 	 \displaystyle \frac{1}{n^p} \sum_{C^{(2)}_{2p} \subseteq C_{2p} } \E[ X_{j_1} X_{j_2} \cdots X_{j_{2p}}] + o(1)  & \text{if }       h=2p, \\
		 	 o(1) & \text{if } h=2p+1.	 
		 	  \end{array}\right. 
	 \end{align} 
Now from the fact (\ref{eqn:opposite_sign_card}), the above equation will be
\begin{align}\label{eqn:limesd_Odd_s} 
\lim_{ \substack{{n\to\infty} \\ {n \mbox{ odd} } }} \E\big[ \tilde{E}_{h}(F_{SC_n}) \big] 
&=  \left\{\begin{array}{ll} 
		 	 \displaystyle \lim_{n \to \infty} \frac{1}{n^p} \Big[\frac{n^p}{2^p} \binom{2p}{p} p! \Big]& \text{if }       h=2p, \\
		 	 0 & \text{if } h=2p+1,
		 	  \end{array}\right. \nonumber \\
& =  \left\{\begin{array}{ll} 
		 	 \frac{(2p)!}{2^p p!} & \text{if }       h=2p, \\
		 	 0 & \text{if } h=2p+1,	 
		 	  \end{array}\right. 
	 \end{align}
 where the first equality arises due to the opposite sign pair matching among the entries of $(j_1, j_2,\ldots, j_{2p})$. 
  Out of $2p$ many entries, $p$ many entries can be freely chosen with a positive sign in $\binom{2p}{p}$ many ways. After free choice of $p$ entries with a positive sign, the rest of $p$ entries with a negative sign can permute among themselves in $p!$ ways. Hence opposite sign pair matching among the entries of $(j_1, j_2, \ldots, j_{2p})$ happen in $\frac{n^p}{2^p} \binom{2p}{p} p!$ ways. 
  Here $1/2^p$ arises because of  $1\leq j_k \leq \lfloor n/2 \rfloor$.

  Now suppose $n$ is even.  Then from  (\ref{trace,formula_sp}), we get
\begin{align*} 
\E\big[ \tilde{E}_{h}(F_{SC_n}) \big] = \frac{1}{2n^{\frac{h}{2} }} \sum_{k=0}^{h}\binom{h}{k} \Big[ \E[Y_k]  \sum_{ C_k} \E[X_{J_{k}}] + \E[\tilde{Y}_k] \sum_{\tilde{C}_k}  \E[X_{J_{k}}] \Big] = \frac{1}{2n^{\frac{h}{2} }} \sum_{k=0}^{h}\binom{h}{k} (T_k +\tilde{T}_k), \mbox{ say},
\end{align*} 
where $ C_{k}$, $\tilde{C}_k$, $Y_k, \tilde{Y}_k, X_{J_k}$ are as defined in  (\ref{def:A_2p_sESD}) for each $k=0,1, \ldots, h$.

First note from (\ref{def:A_2p_sESD}) that $Y_h=2$ and $\tilde{Y}_h=0$. Now from the arguments similar to those used to establish (\ref{eqn:moment_symOdd}), observe that only the terms $T_h$ and $\tilde{T}_h$ can have contribution of the order $O(n^{\lfloor \frac{h}{2} \rfloor})$. But $\tilde{Y}_h=0$, therefore
 \begin{align}\label{eqn:moment_symEven}
 \E\big[ \tilde{E}_{h}(F_{SC_n}) \big] \mathbb I_{\{ n \ even\}}  =
 \left\{\begin{array}{ll} 
		 	 \displaystyle \frac{1}{n^p} \sum_{C^{(2)}_{2p} \subseteq C_{2p} } \E[ X_{j_1} X_{j_2} \cdots X_{j_{2p}}] + o(1),  & \text{if }       h=2p, \\
		 	 o(1) & \text{if } h=2p+1,	 
		 	  \end{array}\right. 
 \end{align}
where $C^{(2)}_{2p}$ is as defined in (\ref{eqn:C^2_k_slsd}). Since the right side of (\ref{eqn:moment_symEven}) is same as the right side of  (\ref{eqn:moment_symOdd}), thus $\lim_{ \substack{{n\to\infty} \\ {n \mbox{ even} } }} \E\big[ \tilde{E}_{h}(F_{SC_n}) \big] $ is equal to the right side of (\ref{eqn:limesd_Odd_s}). Hence  $\lim_{n\to\infty} \E\big[ \tilde{E}_{h}(F_{SC_n}) \big] $ is equal to (\ref{eqn:limesd_Odd_s}).
   
   Observe from the limit of $ \E\big[ \tilde{E}_{h}(F_{SC_n}) \big]$ ((\ref{eqn:limesd_Odd_s})) that the limiting sequence satisfy the condition $(M_2)$ and it corresponds to the standard Gaussian distribution.

Now we check condition $(M_3)$. First note that to prove $(M_3)$, it is sufficient to show that 
\begin{equation} \label{eqn:EO(-2)_sESD}
\E \big[ \tilde{E}_h(F_{SC_n})-\E[ \tilde{E}_h(F_{SC_n})] \big]^4 = O(n^{-2}) \ \mbox{ for all } h \in \mathbb{N}.
\end{equation}
First suppose $n$ is odd. Then from (\ref{trace,formula_sp}), we get
\begin{align} \label{eqn:E^4_odd_sESD}
\E \big[ \tilde{E}_{h} (F_{SC_n})-\E[ \tilde{E}_{h} (F_{SC_n})] \big]^4 
& = \frac{1}{n^{2h} } \E \Big[ \sum_{k=0}^{h}\binom{h}{k} \sum_{C_{k}} \big( X_0^{h-k} X_{J_{k}}  - \E[X_0^{h-k} X_{J_{k}}] \big)   \Big]^4 \nonumber \\
 & \leq \frac{S_1}{n^{2h} } \Big[ \sum_{k=0}^{h} \sum_{C_{k}, C_{k}, C_{k}, C_{k} }  \E \Big\{ \prod_{r=1}^4 \big( X_0^{h-k} X_{J^r_{k}}  - \E[X_0^{h-k} X_{J^r_{k}}] \big) \Big\} \Big],
\end{align}
where $S_1$ is a positive constant (depends only on $h$) 
and
$J^r_{k} = (j^r_1, \ldots, j^r_k) \in C_{k}, \  X_{J^r_k} = X_{j^r_1}  \cdots X_{j^r_k}.$
Note that the inequality in (\ref{eqn:E^4_odd_sESD}) arises due to the fact 
\begin{equation*} 
|x_1 +x_2 + \cdots + x_n|^s \leq 2^{s-1} (|x_1|^s +|x_2|^s +\cdots +|x_n|^s) \ \mbox{ for all } s \in \mathbb{N}.
\end{equation*}
Now following the arguments similar to those used to establish (\ref{eqn:E^4_O(4p-2)_r}), for each $k=0,1, \ldots, h$, we get
\begin{align*} 
 \sum_{C_{k}, C_{k}, C_{k}, C_{k} }  \Big|  \E \Big\{ \prod_{r=1}^4 \big( X_0^{h-k} X_{J^r_{k}}  - \E[X_0^{h-k} X_{J^r_{k}}] \big) \Big\}  \Big|  \leq O(n^{2k-2}).
\end{align*}
Finally, on combining the above inequality with (\ref{eqn:E^4_odd_sESD}), we get
(\ref{eqn:EO(-2)_sESD}) for odd values of $n$.  

The idea of derivation of  (\ref{eqn:EO(-2)_sESD}) for even values of $n$ is similar to the $n$ odd case, so we omit it. This completes the proof of Theorem \ref{thm:symESD}.
\end{proof}
 The following remark provides an alternate trace formula for symmetric circulant matrices, which will be required to prove Theorem \ref{thm:sym_JC}.
 \begin{remark}\label{rem:trace_o(1)_sESD}
Suppose $\{SC_n^{(i)}\}_{1\leq i \leq p}$ are the symmetric circulant matrices with entries $\{\frac{X_j^{(i)} }{\sqrt n};{j\geq 0}\}$. Then using the arguments similar to those used to establish  (\ref{eqn:moment_symOdd}) and (\ref{eqn:moment_symEven}), the trace formula (\ref{trace,formula_sp}) for $\{SC_n^{(i)}\}_{1\leq i \leq p}$ can be written as
\begin{equation}\label{trace,formula_compact_P_sp}
\frac{1}{n}\Tr[ SC_n^{(1)} SC_n^{(2)} \cdots SC_n^{(p)}]
   = \left\{\begin{array}{ccc} 	 
		 	\displaystyle \frac{1}{n^{\frac{p}{2} }} \sum_{C_{p}} X_{j_1}^{(1)} X_{j_2}^{(2)} \cdots X_{j_p}^{(p)} + \theta_p  & \text{if}& \mbox{ $n$ is odd},\\
			\displaystyle \frac{1}{n^{\frac{p}{2} }} \sum_{ C_{p}} X_{j_1}^{(1)} X_{j_2}^{(2)} \cdots X_{j_p}^{(p)} + \gamma_p & \text{if}& \mbox{ $n$ is even}, 	 	 
		 	  \end{array}\right.	
 \end{equation}
where $\E[\theta_p] = \E[\gamma_p]= o(1)$ and $C_p$ is as in (\ref{def:A_2p_sESD}). 
  \end{remark}
\section{Proofs of Theorem \ref{thm:rev_JC} and Theorem \ref{thm:sym_JC}}\label{sec:JC_revsym}
First, we state some definitions and notation which are required for the proof of Theorem \ref{thm:rev_JC}.
\begin{definition} \label{def:half_indep}
Suppose $(\mathcal{A}, \phi)$ is a non-commutative probability space and $\{a_i\}_{i \in I} \subset \mathcal{A}$. Then
	\begin{enumerate}
\item [(i)]We say $\{a_i\}_{i \in I} $ are \textit{half-commutative} if $ {a_ia_ja_k=a_ka_ja_i}$, for all $ i, j, k
\in I$. Moreover, if $\{a_i\}_{i \in I}$ are half-commutative, then $a_i^2$ commutes
with $a_j$ and $a_j^2$ for all $i,j \in I$.
 \item [(ii)] Suppose $a= a_{i_1} a_{i_2} \cdots a_{i_k}$ for any $k \geq1$
and any $\{i_k\} \subset I$. Now for any $ i \in I$, let $E_i(a)$ and $O_i(a)$
denote the number of times $a_{i}$ has occurred in the even
positions and in the odd positions in $a$, respectively. Then $a$ is called \textit{symmetric} with respect to $\{
a_i\}_{i \in I}$ if $E_i(a)=O_i(a)$ for all $i \in I$. Otherwise it is called \textit{non-symmetric}. 
\item [(iii)] Let $\{a_i\}_{i \in I}$ be a family of half commuting elements. Then  $\{a_i\}_{i \in I}$ are called \textit{half independent} if $\{a_i^2\}_{ i \in I}$ are
independent and $\phi(a)=0$ whenever $a$ is non-symmetric with respect to $\{
a_i\}_{i \in I}$.
	\end{enumerate}
\end{definition}
The following example provides a special collection of half independent elements which will appear in the proof of Theorem \ref{thm:rev_JC}.
\begin{example}   \textbf{(Example 2.4, \cite{banica_RC_JC_exam})}.
     \label{exam:example_rJC} 
Let $\big(\mathcal{M}_n, \phi_n=\frac{1}{n} \E\Tr \big)$ be the non-commutative probability space as defined in Definition \ref{def:non-commu} and $\{\eta_i\}$ be a family of independent complex standard Gaussian random variables. Define $a_i \in \big(\mathcal{M}_2, \phi_2\big)$ as
$$ 
a_i= \left(\begin{array}{cccccc}
0 & \eta_i \\
\bar\eta_i & 0
\end{array}\right).
$$
Then $\{a_i\}$ are half independent. Note that each $a_i$ has same distribution as the symmetrized Rayleigh variable and $\{a_i\}$ are half-commutative.
\end{example}
Now let $q_{a}=a_{t_1} a_{t_2} \cdots a_{t_h}$ be a monomial from a collection $\{a_i\}_{1 \leq i \leq m}$. If we suppose $q_{a}$ is symmetric, then $h$ will be even and
$$
q_{a}= \left(\begin{array}{cccccc}
\eta_{t_1} \bar\eta_{t_2} \eta_{t_3} \bar\eta_{t_4} \cdots \eta_{t_{h-1}} \bar\eta_{t_h} & 0 \\
0 & \bar\eta_{t_1} \eta_{t_2} \bar\eta_{t_3} \eta_{t_4} \cdots \bar\eta_{t_{h-1}} \eta_{t_h}
\end{array}\right).
$$
Without loss of generality, assume that for each $i=1,2, \ldots,m$, $a_i$ appears $d_i$ many time in $q_{a}$, where $d_i \geq 0$. Then from the symmetric property of $q_{a}$, $d_i$ will be even, say $2s_i$ and 
$$
q_{a}= \left(\begin{array}{cccccc}
|\eta_{1}|^{2s_1} |\eta_{2}|^{2s_2}  \cdots |\eta_{m}|^{2s_m} & 0 \\
0 & |\eta_{1}|^{2s_1} |\eta_{2}|^{2s_2}  \cdots |\eta_{m}|^{2s_m}
\end{array}\right).
$$ 
Since  $\{\eta_i\}$ is an independent complex standard Gaussian family, we get $\phi[ q_{a}]  = s_1! s_2! \cdots s_m!.$ 

Now suppose $q_{a}$ is non-symmetric.
Then there exists a $i_0 \in \{1,2, \ldots, m\}$ such that $E_{i_0}(q_{a}) \neq O_{i_0}(q_{a})$ and therefore either $\Tr(q_a)=0$ (if $h$ is odd) or a factor  of type $(\eta_{t_{i_0}})^r (\bar\eta_{t_{i_0}})^s$ with $r \neq s$ will appear in each term of $\Tr(q_{a})$. Since $\{\eta_i\}$ is an independent complex standard Gaussian family, $\E[(\eta_{t_{i_0}})^r (\bar\eta_{t_{i_0}})^s]$ will be zero and thus from the independent property of $\{\eta_i\}$, $\phi(q_{a})$ will be zero.

Hence
\begin{align} \label{eqn:phi(amono)_rlsd}
\phi[a_{t_1} a_{t_2} \cdots a_{t_h}] 
  & = \left\{\begin{array}{ll} 	 
		 	s_1! s_2! \cdots s_m! & \text{if }  a_{t_1} a_{t_2} \cdots a_{t_h} \mbox{ is symmetric and for } i=1,2, \ldots, m, \\
		 	 &  \qquad  a_{i} \mbox{ appears $2s_i\geq 0$ many times in } a_{t_1} a_{t_2} \cdots a_{t_h},\\
			0 & \text{if }  a_{t_1} a_{t_2} \cdots a_{t_h} \mbox{ is non-symmetric}. 	 
		 	  \end{array}\right.	
 \end{align}
 
Now we prove Theorem \ref{thm:rev_JC} with the help of the trace formula of reverse circulant matrices.
\begin{proof}[Proof of Theorem \ref{thm:rev_JC}] 
 First note that the reverse circulant matrices half commute. Recall from Definition \ref{def:joint_convergence} that to prove the joint convergence of $\{RC_n^{(i)} \}_{ 1 \leq i \leq m}$, it is suffices to find out the limit of $\phi_n(q)$, where $q= RC^{(t_1)}_n RC^{(t_2)}_n \cdots RC^{(t_h)}_n$ is a monomial in $\{RC_n^{(i)} \}_{ 1 \leq i \leq m}$ and $\phi_n(\cdot) = 1/n \E[\Tr(\cdot)]$.
 Note from the proof of $(M_1)$ of Theorem \ref{thm:revESD} (see (\ref{eqn:ESDo(1)_rESD})) that if $h$ is odd, then $ \phi_n(q)$ will be zero, as $n \to \infty$.  So we consider only the case when $h$ is even, say $2p$. 
From the trace formula (\ref{trace formula RC_n}), we have
\begin{align*} 
 \phi_n(q) =    \frac{1}{n} \E\big[  \Tr\big( RC^{(t_1)}_n RC^{(t_2)}_n \cdots RC^{(t_{2p})}_n \big) \big]
 = \frac{1}{n^{p }}  \sum_{A_{2p}} \E[ X^{(t_1)}_{j^{t_1}_{1}} X^{(t_2)}_{j^{t_2}_{2}} \cdots X^{(t_{2p})}_{j^{t_{2p}}_{2p}} ],
\end{align*}  
 where 
  $J_{2p} =(j^{t_1} _{1}, j^{t_2} _{2}, \ldots, j^{t_{2p}} _{2p} ) \in A_{2p}$ with $A_{2p}$ as in (\ref{def:A_2p_rESD}).  Without loss of generality, assume that $RC_n^{(i)}$ appears $d_i$ many times in $q$, where $d_i \geq 0$ and $\sum_{i=1}^m d_i=2p$.  Furthermore, for $i=1,2, \ldots, m$, suppose $\{j^{i} _{1}, j^{i} _{2}, \ldots, j^{i} _{d_i}\}$ are the entries of $J_{2p}$ which are corresponding to $RC_n^{(i)}$ in $q$. Then  from the independence of $\{RC_n^{(i)}\}$, the above equation will be
\begin{align} \label{eqn:phi_O1_rJC}
 \phi_n(q) 
 = \frac{1}{n^{p }}  \sum_{A_{2p}}  \E \big[  X^{(1)}_{j^{1} _{1}} X^{(1)}_{j^{1} _{2}} \cdots X^{(1)}_{j^{1} _{{2s}_1}} \big]   \E \big[  X^{(2)}_{j^{2} _{1}} X^{(2)}_{j^{2} _{2}} \cdots X^{(2)}_{j^{2} _{{2s}_2}} ]   \cdots  \E \big[ X^{(m)}_{j^{m} _{1}} X^{(m)}_{j^{m} _{2}} \cdots X^{(m)}_{j^{m} _{{2s}_m}} \big] = \frac{1}{n^{p }} Z_p \mbox{ say}.
\end{align}
Now observe from (\ref{eqn:condition_ESD}) and (\ref{eqn:E_bound_rEven}) that a typical term $Z_p$ of (\ref{eqn:phi_O1_rJC}) has the maximum contribution when the following conditions hold: 
 	\begin{enumerate}
\item [(i)] the sets $\{j^{1} _{1}, j^{1} _{2}, \ldots, j^{1} _{d_1}\}, \{j^{2} _{1}, j^{2} _{2}, \ldots, j^{2} _{d_2}\}, \ldots, \{j^{m} _{1}, j^{m} _{2}, \ldots, j^{m} _{d_m} \}$ are pairwise disjoint, 
 \item [(ii)] for each $i=1,2, \ldots, m$,  $d_i$ is even, say $d_i=2s_i$, and
\item [(iii)] for each $i=1,2, \ldots, m$, the elements of the set $\{j^{i} _{1}, j^{i} _{2}, \ldots, j^{i} _{d_i}\}$ are odd-even pair matching (see Definition \ref{def:odd-even_match_rESD}) among themselves in $J_{2p}$.
	\end{enumerate}
Under the above situation, the contribution in $Z_p$  will be of the order $O(n^{s_1+s_2+ \cdots+s_m}) = O(n^{p})$. In other situations, the contribution in $Z_p$  will be the order $o(n^{p})$.	

Now if $q$ is non-symmetric, then there exists a $i_0 \in \{1,2, \ldots, m\}$ such that $E_{i_0}(q) \neq O_{i_0}(q)$. Then the elements of $\{j^{i_0} _{1}, j^{i_0} _{2}, \ldots, j^{i_0} _{d_{i_0}}\}$ can not be odd-even pair matched among themselves in $J_{2p}$.
So the condition (iii) fails and hence this case has contribution of the order $o(1)$ in $\phi_n(q)$.
 
 Now suppose $q$ is symmetric. In this situation, the set $\{j^{i} _{1}, j^{i} _{2}, \ldots, j^{i} _{d_i}\}$ can be odd-even pair matched among themselves for each $i=1,2, \ldots,m$. Hence from the observations (i), (ii) and (iii), for symmetric monomial $q$, (\ref{eqn:phi_O1_rJC}) will be
\begin{align*} 
 \lim_{ n\to \infty} \phi_n(q) 
 & = \lim_{ n\to \infty} \frac{1}{n^p} [ n^{s_1} s_1! n^{s_2} s_2! \cdots n^{s_m} s_m!] 
  = s_1! s_2! \cdots s_m!,
\end{align*}
where the factors $(n^{s_i} s_i!)$'s arise due to the odd-even pair matching among the entries of $(j^{i} _{1}, j^{i} _{2}, \ldots, j^{i} _{2s_i})$.

Suppose $\{a_i\}_{1 \leq i \leq m}$ is a collection of variables as given in Example \ref{exam:example_rJC}. Then from the above equation and (\ref{eqn:phi(amono)_rlsd}), we get that the joint limit of $\{RC_n^{(i)} \}_{ 1 \leq i \leq m}$ is $\{a_i\}_{1 \leq i \leq m}$, which are half independent. This completes the proof of Theorem \ref{thm:rev_JC}.
\end{proof}
\begin{proof}[Proof of Theorem \ref{thm:sym_JC}] 
The idea of the proof is same as the proof of Theorem \ref{thm:rev_JC}. So, it is enough to calculate the limit of $\phi_n(q)$ when $q$ is an even degree (say $2p$) monomial in $\{SC_n^{(i)} \}_{ 1 \leq i \leq m}$. 
 Since symmetric circulant matrices are commutative, $q$ can be written as $q= (SC^{(1)}_n)^{d_1} (SC^{(2)}_n)^{d_2} \cdots (SC^{(m)}_n)^{d_m}$, where $d_i \geq 0$ and $\sum_{i=1}^m d_i=2p$.
Now from (\ref{trace,formula_compact_P_sp}), for both odd and even values of $n$, we have
\begin{align} \label{eqn:E,q_even_s}
 \phi_n(q)     =   \frac{1}{n^p}  \sum_{ C_{2p}} \E[ X^{(1)}_{J_{d_1}} X^{(2)}_{J_{d_2}} \cdots X^{(m)}_{J_{d_m}}] +o(1)
   = \frac{1}{n^p}  \sum_{C_{2p}} \E[ X^{(t_1)}_{J_{m_1}}] \E[ X^{(t_2)}_{J_{m_2}}] \cdots \E[X^{(t_r)}_{J_{m_r}}] +o(1),
\end{align}
where $C_{2p}$ is as defined in (\ref{def:A_2p_sESD}) and for each $i=1, 2, \ldots, m$,
\begin{align*} 
J_{d_i} & = (j^{i} _{1}, \ldots, j^{i} _{d_i}), \  X^{(i)}_{J_{d_i}}= X^{(i)}_{j^{i} _{1}} \cdots X^{(i)}_{j^{i} _{d_i}} \mbox{ and }
J_{2p} = (j^{1} _{1}, \ldots, j^{1} _{d_1}, j^{2} _{1}, \ldots, j^{2} _{d_2}, \ldots, j^{m} _{1}, \ldots, j^{m} _{d_m}) \in C_{2p}. 
\end{align*}
Note that the second equality in (\ref{eqn:E,q_even_s}) arises due to the independence property of $\{SC_n^{(i)} \}_{ 1 \leq i \leq m}$.
Now observe from the fact (\ref{eqn:opposite_sign_card}) that a typical term of (\ref{eqn:E,q_even_s}) has the maximum contribution when the following conditions hold: 
\begin{enumerate}
	\item[(i)] the sets $\{j^{1} _{1}, j^{1} _{2}, \ldots, j^{1} _{d_1}\}, \{j^{2} _{1}, j^{2} _{2}, \ldots, j^{2} _{d_2}\}, \ldots, \{j^{m} _{1}, j^{m} _{2}, \ldots, j^{m} _{d_m} \}$ are pairwise disjoint, and
\item[(ii)] for each $i=1,2, \ldots, m$, $d_i$ is even, say $d_i=2s_i$ and $J_{d_i}$ is  opposite sign pair matched (see Definition \ref{def:opposite_sign}).
	\end{enumerate}  
In this case, the contribution will be of the order $O(n^{s_1+s_2+ \cdots+s_m}) = O(n^{p})$. In other situations, the contribution will be the order $o(n^{p})$. Hence from (\ref{eqn:E,q_even_s}), we get
\begin{align} \label{eqn:limE,q_sJC}
 \lim_{ n\to \infty} \phi_n(q) 
 & = \lim_{ n\to \infty} \frac{1}{n^p} \big[ n^{s_1} N_{s_1} n^{s_2} N_{s_2} \cdots n^{s_m} N_{s_m} \big]
 = N_{s_1} N_{s_2} \cdots N_{s_m},
\end{align}
where $N_{s_i}= \frac{(2s_i)!}{2^{s_i} s_i!}$.
Note that the factors $(n^{s_i}N_{s_i})$'s arise in the above equation due to the opposite sign pair matching among the entries of $(j^{i} _{1}, j^{i} _{2}, \ldots, j^{i} _{2s_i})$.

Now, suppose $\{b_i\}_{1 \leq i \leq m}$ is a collection of independent standard Gaussian variables, then
\begin{align*}
\E\big( (b_{1})^{d_1} (b_{2})^{d_2} \cdots (b_{m})^{d_m} \big) 
 &  = \left\{\begin{array}{ll} 	 
		 	N_{s_1} N_{s_2} \cdots N_{s_m} & \text{if } d_i=2s_i\geq 0, \ \forall \ i = 1, 2, \ldots, m, \\
			0 & \text{otherwise}. 	 	 
		 	  \end{array}\right.	
\end{align*} 
Hence it follows from  the above expression and (\ref{eqn:limE,q_sJC}) that the joint limit of $\{SC_n^{(i)} \}_{ 1 \leq i \leq m}$ is $\{b_i\}_{1 \leq i \leq m}$, which are independent. This completes the proof of Theorem \ref{thm:rev_JC}.
\end{proof}

\noindent\textbf{Acknowledgement:} 
 The author would like to express his heartfelt gratitude to Prof. Koushik Saha for reading the draft and providing valuable suggestions. The work of author is partially supported by UGC Doctoral Fellowship, India.


\providecommand{\bysame}{\leavevmode\hbox to3em{\hrulefill}\thinspace}
\providecommand{\MR}{\relax\ifhmode\unskip\space\fi MR }
\providecommand{\MRhref}[2]{%
  \href{http://www.ams.org/mathscinet-getitem?mr=#1}{#2}
}
\providecommand{\href}[2]{#2}

\end{document}